# Bahadur–Kiefer theory for sample quantiles of weakly dependent linear processes

RAFAŁ KULIK

[1] *University of Sydney, School of Mathematics and Statistics, F07, University of Sydney, NSW 2006, Australia. E-mail: rkuli@maths.usyd.edu.au*

In this paper, we establish the Bahadur–Kiefer representation for sample quantiles for a class of weakly dependent linear processes. The rate of approximation is the same as for i.i.d. sequences and is thus optimal.

*Keywords:* Bahadur representation; empirical processes; law of the iterated logarithm; linear processes; sample quantiles; strong approximation

## 1. Introduction

Consider the class of stationary linear processes

$$X_i = \sum_{k=0}^{\infty} c_k \epsilon_{i-k},$$

where $\{\epsilon_i, i \in \mathbb{Z}\}$ is an i.i.d. sequence and $\sum_{k=0}^{\infty} |c_k| < \infty$. Assume that $X_1$ has continuous distribution function $F(x) = P(X_1 \le x)$ and let $f$ and $Q$ denote the associated density and quantile function, respectively. Given a sample $X_1, \ldots, X_n$, let $F_n(x) = \frac{1}{n} \sum_{i=1}^{n} 1_{\{X_i \le x\}}$ and let $Q_n(y)$ denote the corresponding empirical quantile function. Define

$$\beta_n(x) = n^{1/2}(F_n(x) - F(x)), \qquad x \in \mathbb{R},$$
$$q_n(y) = n^{1/2}(Q(y) - Q_n(y)), \qquad y \in (0,1),$$

the *general empirical* and the *general quantile* process, respectively. With $U_i = F(X_i)$, $i \ge 1$, let $E_n(x) = \frac{1}{n} \sum_{i=1}^{n} 1_{\{U_i \le x\}} = F_n(Q(x))$ and $U_n(y)$ be the uniform empirical distribution and the uniform empirical quantile function, respectively. Let

$$\alpha_n(x) = n^{1/2}(E_n(x) - x), \qquad x \in (0,1),$$
$$u_n(y) = n^{1/2}(y - U_n(y)), \qquad y \in (0,1),$$

be the corresponding *uniform empirical* and *uniform quantile* process, respectively.

Assume, initially, that $X_i$, $i \geq 1$, are i.i.d. Fix $y \in (0,1)$. Let $I_y$ be a neighborhood of $Q(y)$. Assuming that $\inf_{x \in I_y} f(x) > 0$ and $\sup_{x \in I_y} |f^{'}(x)| < \infty$, Bahadur [2] obtained the representation

$$f(Q(y))q_n(y) - \alpha_n(y) =: R_n(y),\tag{1}$$

where $R_n(y) = O_{a.s}(n^{-1/4}(\log n)^{1/2}(\log\log n)^{1/4})$. Later, Kiefer [11] proved that this can be strengthened to

$$R_n(y) = O_{a.s}(n^{-1/4}(\log\log n)^{3/4}),\tag{2}$$

which is the optimal rate. Continuing his study, Kiefer [12] established the uniform version of (1), later referred to as the *Bahadur–Kiefer representation*. For $0 \leq a < b \leq 1$,

$$\sup_{y \in (a,b)} |f(Q(y))q_n(y) - \alpha_n(y)| =: R_n,\tag{3}$$

where

$$R_n = O_{a.s}(n^{-1/4}(\log n)^{1/2}(\log\log n)^{1/4}).\tag{4}$$

The above rate is also optimal. Kiefer obtained his result assuming

(K1)  $\sup_{y \in (a,b)} |f^{'}(Q(y))| < \infty$;
(K2)  $\inf_{y \in (a,b)} f(Q(y)) > 0$.

We shall refer to (K1)–(K2) as the *Kiefer conditions*. In particular, if $a = 0$ and $b = 1$, (K1)–(K2) imply that $f$ has a finite support.

Further, Csörgő and Révész [5] relaxed the conditions for Kiefer's result (3) and introduced the *Csörgő–Révész conditions* (see Section 2; also cf. [4], Theorem 3.2.1).

The main purpose of this paper is to obtain Bahadur–Kiefer-type representations for sample quantiles of linear processes with the *optimal* rate. Some results are available for weakly dependent random variables. For $\phi$-mixing sequences, under (K1)–(K2), the optimal rates have been obtained in [1]. These results were improved in [7] and [18] through less restrictive mixing rates and Csörgő–Révész-type conditions. The rate of approximation was $R_n = O_{a.s}((\log n)^{-\lambda})$ for some $\lambda > 0$.

Mixing is rather hard to verify and requires additional assumptions. In particular, to obtain a strong mixing for linear processes, both regularity of the density of $\epsilon_1$ and some constraints on the $c_k$'s are required (see, e.g., [6] or [15]). Nevertheless, even if we are able to establish strong mixing, we do not attain the optimal rate in the Bahadur–Kiefer representation.

Another way of looking at linear processes is to approximate the sequence $\{X_i, i \geq 1\}$ by a sequence with finite memory and then to use the classical Bernstein blocking technique. Hesse [8] obtained the Bahadur representation (1) with rate $R_n(y) = O_{a.s}(n^{-1/4+\lambda})$ for some $\lambda > 0$ via this technique. The method avoids some assumptions on the density of $\epsilon_1$, but it leads to restrictive constraints on the $c_k$'s and does not lead to the optimal rates.

Both the blocking technique and mixing require some strong assumptions. To overcome such restrictions, Ho, Hsing, Mielniczuk and Wu (see [9, 10, 16, 17]) developed martingale-based methods. In particular, for a class of linear processes, Wu [16] obtained the exact rate (2) in the Bahadur representation (1). He also studied the Bahadur–Kiefer representation, but was not able to obtain the optimal rate (4) due to the lack of an appropriate version of the law of the iterated logarithm for empirical processes.

In this paper, we shall combine Bernstein's blocking technique with Wu's method to obtain the optimal rate in the Bahadur–Kiefer representation (3) under quite mild conditions on the $c_k$'s. The result is stated in Theorem 2.1. The methodology involves the recent strong approximation result of Berkes and Horváth [3]. Since part of our computation follows their proof, we include it in the Appendix. Further, we shall obtain the optimal rate, under general conditions on $F$, on the whole interval $(0, 1)$ by considering an appropriately weighted process $f(Q(y))q_n(y) - \alpha_n(y)$, as in Theorem 2.2.

Throughout the paper, $C$ will denote a generic constant which may be different at each appearance. Also, we write $a_n \sim b_n$ if $\lim_{n\to\infty} a_n/b_n = 1$. For any stationary sequence $\{Z_i, i \geq 1\}$ of random variables, $Z$ will be a random variables with the same distribution as $Z_1$.

Some further notation is required. Let

$$b_n = n^{-1/4}(\log n)^{1/2}(\log\log n)^{1/4}$$

and

$$\lambda_n = n^{-1/2}(2\log\log n)^{1/2}.$$

For any function $h(x)$ defined on $\mathbb{R}$ and $x < y$, we write $h(x, y) := h(y) - h(x)$.

## 2. Results

Assume that the following moment and dependence conditions hold. For $\alpha \geq 2$,

$$\mathrm{E}|\epsilon|^\alpha < \infty \tag{5}$$

and for some $\rho \in (0, \frac{1}{2})$,

$$\sum_{k=i}^\infty c_k^2 = O(i^{-2/\rho}(\log i)^{-3}). \tag{6}$$

**Theorem 2.1.** *Assume* (5), (6) *and* (K1)–(K2). *Furthermore, assume that for* $f_\epsilon$, *a density of* $\epsilon$, *we have*

$$\sup_{x\in\mathbb{R}}(f_\epsilon(x) + |f_\epsilon'(x)| + |f_\epsilon''(x)|) < \infty. \tag{7}$$

*Then,*

$$\sup_{y\in(a,b)} |f(Q(y))q_n(y) - \alpha_n(y)| = O_{a.s}(n^{-1/4}(\log n)^{1/2}(\log\log n)^{1/4}).$$

To obtain the bound without the Kiefer conditions (K1)–(K2), we shall consider the following *Csörgő–Révész conditions*:

(CsR1) $f'$ exists on $(a,b)$, where $a = \sup\{x: F(x) = 0\}$ and $b = \inf\{x: F(x) = 1\}$, $-\infty \le a < b \le \infty$;
(CsR2) $\inf_{x\in(a,b)} f(x) > 0$;
(CsR3) (i) $f(Q(y)) \sim y^{\gamma_1}L_1(y^{-1})$ as $y \downarrow 0$;
(CsR3) (ii) $f(Q(y)) \sim (1-y)^{\gamma_2}L_2((1-y)^{-1})$ as $y \uparrow 1$;
(CsR4) (i) $0 < A := \lim_{y\downarrow 0} f(Q(y)) < \infty$, $0 < B := \lim_{y\uparrow 1} f(Q(y)) < \infty$, or
      (ii) if $A = 0$ (resp. $B = 0$), then $f$ is non-decreasing (resp. non-increasing) on an interval to the right of $Q(0+)$ (resp. to the left of $Q(1-)$).

**Theorem 2.2.** *Assume Csörgő–Révész conditions with $\gamma := \min\{\gamma_1, \gamma_2\} \ge 1$. As in Theorem 2.1, assume (5), (6) and (7). Then, for arbitrary $\nu > \max\{2\gamma, 3\gamma - 2\}$,*

$$\sup_{y\in(0,1)} (y(1-y))^{\nu} |f(Q(y))q_n(y) - \alpha_n(y)| = O_{a.s}(n^{-1/4}(\log n)^{1/2}(\log\log n)^{1/4}).$$

# 3. Proof of Theorem 2.1

Let $\mathcal{F}_i = \sigma(\epsilon_i, \epsilon_{i-1}, \ldots)$. Note that $F(x) = \mathrm{E}F_\epsilon(x - X_{i,i-1}) = \mathrm{E}F_\epsilon(x - X_{1,0})$. Also, let $Y_i(x) = F_\epsilon(x - X_{i,i-1}) - F(x)$. We then have

$$\frac{1}{n}\sum_{i=1}^{n}(1_{\{X_i \le x\}} - F(x))$$

$$= \frac{1}{n}\sum_{i=1}^{n}(1_{\{X_i \le x\}} - \mathrm{E}(1_{\{X_i \le x\}}|\mathcal{F}_{i-1})) + \frac{1}{n}\sum_{i=1}^{n}Y_i(x)$$

$$=: M_n(x) + N_n(x).$$

Then $nM_n(x)$, $n \ge 1$, is a martingale and $N_n(x)$ is differentiable.

The plan of our proof is as follows. First, we obtain a strong approximation of the differentiable part $N_n$ by an appropriate Gaussian process. In order to do this, we will replace the original sequence $\{X_i, i \ge 1\}$ with one having finite memory. From that approximation, we will establish the uniform law of the iterated logarithm (ULIL) for the differentiable part $N_n$, which, together with the ULIL for the martingale part $M_n$, will imply the ULIL for the empirical process $\beta_n$ (see Section 3.1).

Then, using a modification of Lemma 13 from [16], we will be able to control increments of the empirical processes $\beta_n$ and $\alpha_n$ (Section 3.2), which, together with the ULIL, will imply the result (Section 3.3).

## 3.1. Approximation of the differentiable part by a Gaussian process and laws of the iterated logarithm

**Proposition 3.1.** *Assume* *(5)*, *(6)* *and* *(7)*. *There then exists a centered Gaussian process* $K(x,t)$ *with* $\mathbb{E}K(x,t)K(y,t') = t \wedge t'\Gamma(x,y)$ *such that*

$$\sup_{0 \le t \le 1} \sup_{x \in \mathbb{R}} |[nt]N_{[nt]}(x) - K(x,[nt])| = o(n^{1/2}(\log n)^{-\lambda}) \qquad \text{almost surely}$$

*for some* $\lambda > 0$.

The law of the iterated logarithm follows from Proposition 3.1 and the ULIL for $K$.

**Corollary 3.2.** *If* *(5)*, *(6)*, *(7)* *and Kiefer conditions are fulfilled, then*

$$\limsup_{n \to \infty} \frac{1}{(2\log\log n)^{1/2}} \sup_{y \in (0,1)} |q_n(y)| = C \qquad \text{almost surely.} \tag{8}$$

**Proof.** Proposition 3.1, together with the ULIL for the martingale $M_n$ (see, e.g., [16], Lemma 7), yields

$$\limsup_{n \to \infty} \frac{1}{(2\log\log n)^{1/2}} \sup_{-\infty < x < \infty} |\beta_n(x)| = C \qquad \text{almost surely.}$$

Consequently,

$$\limsup_{n \to \infty} \frac{1}{(2\log\log n)^{1/2}} \sup_{y \in (0,1)} |u_n(y)|$$
$$= \limsup_{n \to \infty} \frac{1}{(2\log\log n)^{1/2}} \sup_{x \in (0,1)} |\alpha_n(x)| = C \tag{9}$$

almost surely.

With $\Delta_{n,y} = Q_n(y) - Q(y)$, we have

$$\limsup_{n \to \infty} \sup_{y \in (a,b)} \frac{n^{1/2}\Delta_{n,y}}{(\log\log n)^{1/2}}$$

$$\le \limsup_{n \to \infty} \sup_{1 \le k \le n} \frac{n^{1/2}|Q(F(X_{k:n})) - Q(k/n)|}{(\log\log n)^{1/2}}$$

$$+ \limsup_{n \to \infty} \sup_{1 \le k \le n} \sup_{y \in ((k-1)/n, k/n]} \frac{n^{1/2}|Q(k/n) - Q(y)|}{(\log\log n)^{1/2}}.$$

Now, if the Kiefer conditions hold, then $\sup_{y \in (a,b)} |Q''(y)| < \infty$. Using Taylor's expansion, we obtain $Q(\frac{k}{n}) = Q(y) + Q'(y)(\frac{k}{n} - y) + O(\frac{1}{n})$ and

$$Q(F(X_{k:n})) = Q(U_{k:n}) = Q\left(\frac{k}{n}\right) + Q'\left(\frac{k}{n}\right)\left(U_{k:n} - \frac{k}{n}\right) + O_{a.s}\left(\left(U_{k:n} - \frac{k}{n}\right)^2\right).$$

Thus, by (9),

$$\limsup_{n \to \infty} \sup_{y \in (a,b)} \frac{n^{1/2} \Delta_{n,y}}{(\log \log n)^{1/2}} \leq \sup_{y \in (a,b)} |Q'(y)| \limsup_{n \to \infty} \sup_{1 \leq k \leq n} \frac{n^{1/2} |U_{k:n} - k/n|}{(\log \log n)^{1/2}} \leq C. \quad \square$$

Before proving Proposition 3.1, we need some notation. With $\rho$ from (6), define

$$\hat{X}_i := \hat{X}_i(\rho) = \sum_{k=0}^{i^\rho - 1} c_k \epsilon_{i-k}.$$

Without loss of generality, we may assume that $c_0 = 1$. Let $X_{i,i-1} = \sum_{k=1}^{\infty} c_k \epsilon_{i-k}$ and define its truncated version $\hat{X}_{i,i-1} = \sum_{k=1}^{i^\rho - 1} c_k \epsilon_{i-k}$.

By Rosenthal's inequality, for any $\alpha \geq 2$,

$$\mathrm{E}|X_i - \hat{X}_i|^\alpha = \mathrm{E}\left|\sum_{k=i^\rho}^{\infty} c_k \epsilon_{i-k}\right|^\alpha$$

$$\leq C\left(\sum_{k=i^\rho}^{\infty} c_k^2\right)^{\alpha/2} + C \sum_{k=i^\rho}^{\infty} |c_k|^\alpha. \tag{10}$$

This estimate is also true for $\mathrm{E}|X_{i,i-1} - \hat{X}_{i,i-1}|^\alpha$. Thus, replacing $i$ with $i^\rho$ in (6), we have

$$\left\|\sum_{i=1}^{n} X_i - \sum_{i=1}^{n} \hat{X}_i\right\|_\alpha \leq \sum_{i=1}^{\infty} O(i^{-1}(\log i)^{-3/2}) < \infty \tag{11}$$

and the same estimate is valid for $\|\sum_{i=1}^{n} X_{i,i-1} - \sum_{i=1}^{n} \hat{X}_{i,i-1}\|_\alpha$.

Define

$$F_n^*(x) = \frac{1}{n} \sum_{i=1}^{n} F_\epsilon(x - X_{i,i-1}),$$

a conditional empirical distribution function, and

$$\hat{F}_n^*(x) = \frac{1}{n} \sum_{i=1}^{n} F_\epsilon(x - \hat{X}_{i,i-1}),$$

its corresponding version based on the truncated random variables.

Let $\hat{F}_i(x) := P(\hat{X}_i \leq x) = \mathrm{E} f_\epsilon(x - \hat{X}_{i,i-1}), i \geq 1$. Since $f_\epsilon$ exists, we may define $f_n^*(x) = \mathrm{d} F_n^*(x)/\mathrm{d}x$, $\hat{f}_n^*(x) = \mathrm{d} \hat{F}_n^*(x)/\mathrm{d}x$, $\hat{f}_n(x) = \frac{1}{n} \sum_{i=1}^n \mathrm{d}\hat{F}_i(x)/\mathrm{d}x$.

Further, define $\hat{Y}_i(x) = F_\epsilon(x - \hat{X}_{i,i-1}) - \hat{F}_i(x)$ and

$$\hat{N}_n(x) = \frac{1}{n} \sum_{i=1}^n (F_\epsilon(x - \hat{X}_{i,i-1}) - \hat{F}_i(x)) = \frac{1}{n} \sum_{i=1}^n \hat{Y}_i(x).$$

**Lemma 3.3.** *Assume [(5)](#), [(6)](#) and [(7)](#). Then, for some $D_0 > 0$,*

$$\|n\hat{N}_n(x,y)\|_2^2 \leq D_0 n(y-x)^2. \tag{12}$$

**Proof.** We have, for any $-\infty < x < y < \infty$,

$$\begin{aligned}
\|nN_n(x,y)\|_2^2 &= n^2 \mathrm{E} \left| \int_x^y (f_n^*(u) - f(u)) \, \mathrm{d}u \right|^2 \\
&\leq n^2 (y-x)^2 \mathrm{E} \sup_{x \leq u \leq y} |f_n^*(u) - f(u)|^2 \\
&\leq Cn(y-x)^2,
\end{aligned} \tag{13}$$

by Lemma 9 of [[16]](#). Next,

$$\begin{aligned}
\mathrm{E} \sup_{x \leq u \leq y} |\hat{f}_n^*(u) - \hat{f}_n(u)|^2 \leq 3 \Bigg( &\mathrm{E} \sup_{x \leq u \leq y} |\hat{f}_n^*(u) - f_n^*(u)|^2 \\
&+ \mathrm{E} \sup_{x \leq u \leq y} |f_n^*(u) - f(u)|^2 + \sup_{x \leq u \leq y} |f(u) - \hat{f}_n(u)|^2 \Bigg)
\end{aligned} \tag{14}$$

The second term is $O(n^{-1})$, as in [(13)](#). As for the first term, we have

$$\begin{aligned}
\mathrm{E} \sup_{x \leq u \leq y} |\hat{f}_n^*(u) - f_n^*(u)|^2 &= \frac{1}{n^2} \mathrm{E} \sup_{x \leq u \leq y} \left| \sum_{i=1}^n (f_\epsilon(u - X_{i,i-1}) - f_\epsilon(u - \hat{X}_{i,i-1})) \right|^2 \\
&= \frac{1}{n^2} \mathrm{E} \sup_{x \leq u \leq y} \left| \sum_{i=1}^n \int_{u - X_{i,i-1}}^{u - \hat{X}_{i,i-1}} f_\epsilon'(v) \, \mathrm{d}v \right|^2 \\
&\leq \frac{1}{n^2} \sup_{x \in \mathbb{R}} |f_\epsilon'(x)|^2 \mathrm{E} \left( \sum_{i=1}^n |X_{i,i-1} - \hat{X}_{i,i-1}| \right)^2 \\
&\leq \frac{1}{n} \sup_{x \in \mathbb{R}} |f_\epsilon'(x)|^2 \left( \sum_{i=1}^n \|X_{i,i-1} - \hat{X}_{i,i-1}\|_\alpha \right)^2 \leq C \frac{1}{n}
\end{aligned} \tag{15}$$

due to the fact that $\alpha \geq 2$ and the comment following (11). Also,

$$\sup_{x \leq u \leq y} |f(u) - \hat{f}_n(u)|^2 \leq \sup_{x \leq u \leq y} \mathrm{E}|\hat{f}_n^*(u) - f_n^*(u)|^2 \leq C\frac{1}{n}. \qquad (16)$$

Combining (14), (15) and (16), we obtain (12). $\qquad \square$

Next, we derive an exponential inequality for $nN_n(x,y)$.

**Lemma 3.4.** *Assume (5), (6) and (7). Then, for any $z > 0$,*

$$P\left(\left|\sum_{i=1}^n Y_i(x,y)\right| > z\right) \leq C_1 z^{-\alpha} + C_2 \exp(-C_3 z^2/(n(y-x)^2)) + C_4 \exp(-C_5 z/n^\rho),$$

*where $C_1, C_2, C_3, C_4, C_5$ are positive constants.*

**Proof.** Making use of differentiability and the comment following (10), we obtain

$$\|F_\epsilon(x - X_{i,i-1}) - F_\epsilon(x - \hat{X}_{i,i-1})\|_\alpha = O\left(\left(\left(\sum_{k=i^\rho}^\infty c_k^2\right)^{\alpha/2} + \sum_{k=i^\rho}^\infty |c_k|^\alpha\right)^{1/\alpha}\right).$$

Further,

$$|F(x) - \hat{F}_i(x)|^\alpha \leq \mathrm{E}|F_\epsilon(x - X_{i,i-1}) - F_\epsilon(x - \hat{X}_{i,i-1})|^\alpha.$$

Consequently,

$$\|Y_i(x) - \hat{Y}_i(x)\|_\alpha = O\left(\left(\left(\sum_{k=i^\rho}^\infty c_k^2\right)^{\alpha/2} + \sum_{k=i^\rho}^\infty |c_k|^\alpha\right)^{1/\alpha}\right) \qquad (17)$$

and

$$\|nN_n(x) - n\hat{N}_n(x)\|_\alpha \leq C_0 := \sum_{i=1}^\infty i^{-1}(\log i)^{-3/2}.$$

From this and the Markov inequality, we get

$$P(|nN_n(x,y)| > z) \leq C_0 z^{-\alpha} + P(|n\hat{N}_n(x,y)| > z/2).$$

To obtain the bound for the second part, divide $[1,n]$ into blocks $I_1, J_1, I_2, J_2, \ldots, I_M, J_M$ with the same length $n^\rho$, where $\rho$ is defined in (6). Thus, $M \sim n^{1-\rho}$. Let $\hat{U}_k = \sum_{i \in I_k} \hat{Y}_i(x,y)$, $\hat{V}_k = \sum_{i \in J_k} \hat{Y}_i(x,y)$ and $n\hat{N}_n^{(1)} = \sum_{k=1}^M \hat{U}_k$, $n\hat{N}_n^{(2)} = \sum_{k=1}^M \hat{V}_k$. Both $(\hat{U}_1, \ldots, \hat{U}_M)$ and $(\hat{V}_1, \ldots, \hat{V}_M)$ are vectors of independent random variables. Also, $\max_{k=1,\ldots,M} \hat{U}_k \leq [I_k] \leq Cn^\rho$. Equation (12) yields $\|U_k\|_2^2 \leq D_0 n^\rho (y-x)^2$. Recall that $\rho \in (0, \frac{1}{2})$. Applying the result from [13] page 293, to the centered sequence $\hat{U}_1, \ldots, \hat{U}_M$ with $M_n = n^\rho$ and $B_n = D_0 n(y - x)^2$, we obtain

$$P(|n\hat{N}_n^{(1)}| > z) \leq \exp(-z^2/(4D_0 n(y - x)^2)) + \exp(-z/(4n^\rho)).$$

The same applies to $P(|n\hat{N}_n^{(2)}| > z)$ and hence the result follows. □

To state an approximation result, first note that $Y_0(x)$ and $\hat{Y}_i(x)$, $i \geq 1$, are independent. Thus, we obtain

$$|\mathbb{E}Y_0(x)Y_i(y)| \leq |\mathbb{E}Y_0(x)Y_i(y) - Y_0(x)\hat{Y}_i(y)| + |\mathbb{E}Y_0(x)\hat{Y}_i(y)|$$
$$\leq \|Y_i(y) - \hat{Y}_i(y)\|_\alpha + 0.$$

Then, in view of (6) and (17),

$$\Gamma(x, y) = \mathbb{E}Y_0(x)Y_0(y) + \sum_{i=1}^{\infty} (\mathbb{E}Y_0(x)Y_i(y) + \mathbb{E}Y_0(y)Y_i(x)) \qquad (18)$$

is absolutely convergent for all $x, y \in \mathbb{R}$.

Having Lemma 3.4, (18) and $\alpha \geq 2$, we may proceed exactly as in [3] to obtain the approximation result of Proposition 3.1 (see the Appendix for details).

## 3.2. Controlling increments

Let $l_q(n) = (\log n)^{1/q} (\log \log n)^{2/q}$, $q > 2$, and recall that $\sum_{k=0}^{\infty} |c_k| < \infty$.

First, we generalize Lemma 13 of [16], to the real line.

**Lemma 3.5.** *Assume* (5), (7) *and* (K1)–(K2). *Let* $-\infty \leq a < b \leq \infty$. *Then, under Kiefer conditions, for any positive bounded sequence* $d_n$ *such that* $\log n = o(nd_n)$, *we have*

$$\sup_{|x-y| \leq d_n, x, y \in [a,b]} |F_n(x) - F(x) - (F_n(y) - F(y))| = O_{a.s}\left(\frac{\sqrt{d_n \log n}}{\sqrt{n}} + \frac{d_n l_q(n)}{\sqrt{n}}\right). \qquad (19)$$

**Proof.** Define $U_i = F(X_i)$. Clearly, in order to show (19), it suffices to prove that

$$\sup_{|u-v| \leq d_n, u, v \in (a_1, b_1)} |E_n(v) - v - (E_n(u) - u)| = O_{a.s}\left(\frac{\sqrt{d_n \log n}}{\sqrt{n}} + \frac{d_n l_q(n)}{\sqrt{n}}\right). \qquad (20)$$

Let $y \in (a_1, b_1)$. Decompose

$$\frac{1}{n} \sum_{i=1}^{n} (1_{\{U_i \leq y\}} - y)$$

$$= \frac{1}{n}\sum_{i=1}^{n}(1_{\{U_i \leq y\}} - \mathrm{E}(1_{\{U_i \leq y\}}|\mathcal{F}_{i-1}))$$

$$+ \frac{1}{n}\sum_{i=1}^{n}(\mathrm{E}(1_{\{U_i \leq y\}}|\mathcal{F}_{i-1}) - y)$$

$$=: \tilde{M}_n(x) + \tilde{N}_n(x),$$

where

$$\tilde{N}_n(y) = \frac{1}{n}\sum_{i=1}^{n}(\mathrm{E}(1_{\{U_i \leq y\}}|\mathcal{F}_{i-1}) - y) = \frac{1}{n}\sum_{i=1}^{n}(F_\epsilon(Q(y) - X_{i,i-1}) - y)$$

so that

$$\tilde{N}_n'(y) = \frac{1}{n}\sum_{i=1}^{n}(Q'(y)f_\epsilon(Q(y) - X_{i,i-1}) - 1) \tag{21}$$

and

$$\tilde{N}_n''(y) = \frac{1}{n}\sum_{i=1}^{n}(Q''(y)f_\epsilon(Q(y) - X_{i,i-1}))$$

$$+ (Q'(y))^2 \frac{1}{n}\sum_{i=1}^{n}f_\epsilon'(Q(y) - X_{i,i-1}). \tag{22}$$

Define the projection operator $\mathcal{P}_k\xi = \mathrm{E}(\xi|\mathcal{F}_k) - \mathrm{E}(\xi|\mathcal{F}_{k-1})$. Let $g_y(\mathcal{F}_i) = f_\epsilon(Q(y) - X_{i,i-1})$. As in the proof of Lemma 3 of [16], with $\alpha = q$,

$$\|\mathcal{P}_0 g_y(\mathcal{F}_i)\|_\alpha = O(|c_i|),$$

uniformly in $y \in (a_1, b_1)$. The same holds for $g_y(\mathcal{F}_i) = f_\epsilon'(Q(y) - X_{i,i-1})$. Further, under Kiefer conditions, $\sup_{y \in (a_1, b_1)}(|Q'(y)| + |Q''(y)|) < \infty$. Thus, applying Proposition 1 of [16], we have $\max_{y \in (a_1, b_1)}(\|\tilde{N}_n'(y)\|_\alpha + \|\tilde{N}_n''(y)\|_\alpha) = O(n^{-1/2})$. Consequently,

$$\mathrm{E}\left[\max_{y \in (a_1, b_1)} |\tilde{N}_n(y)|^\alpha\right] \leq \mathrm{E}\left[\int_0^1 |\tilde{N}_n'(y)| \, \mathrm{d}y\right]^\alpha = O(n^{-\alpha/2}).$$

Similarly, $\mathrm{E}[\max_{y \in (a_1, b_1)} |\tilde{N}_n'(y)|^\alpha] = O(n^{-\alpha/2})$. Thus, by the same argument as in Lemma 9 of [16] we have $\sup_{y \in (a_1, b_1)} |\tilde{N}_n'(y)| = o_{a.s}(l_q(n)/\sqrt{n})$. Further,

$$\sup_{|u-v| \leq d_n, u, v \in (a_1, b_1)} |\tilde{N}_n(v) - \tilde{N}_n(u)| \leq d_n \sup_{y \in (a_1, b_1)} |\tilde{N}_n'(y)| = o_{a.s}(d_n l_q(n)/\sqrt{n}).$$

This, together with appropriate estimates for the martingale part, yields (20).  $\square$

### 3.3. Conclusion of the proof of Theorem 2.1

We apply Lemma 3.5 with $d_n = \lambda_n$. The second part in (19) is then negligible. Thus, we have

$$\sup_{x \in \mathbb{R}} \sup_{|x-y| \leq \lambda_n} |\beta_n(x) - \beta_n(y)| = O_{a.s}(b_n).$$

On account of (8), it yields

$$\sup_{y \in (0,1)} |\beta_n(Q_n(y)) - \beta_n(Q(y))| = O_{a.s}(b_n).$$

Equivalently,

$$n^{1/2} \sup_{y \in (0,1)} |F_n(Q_n(y)) - F(Q_n(y)) - (F_n(Q(y)) - F(Q(y)))| = O_{a.s}(b_n).$$

Since $|F_n(Q_n(y)) - F(Q(y))| \leq 1/n$, we obtain

$$n^{1/2} \sup_{y \in (0,1)} |F(Q_n(y)) - F(Q(y)) - (F(Q(y)) - F_n(Q(y)))| = O_{a.s}(b_n).$$

Set $\Delta_{n,y} = Q_n(y) - Q(y)$. Using the Taylor expansion $F(Q_n(y)) = F(Q(y)) + f(Q(y))\Delta_{n,y} + O_{a.s}(\Delta_{n,y}^2)$, we complete the proof of Theorem 2.1. □

## 4. Proof of Theorem 2.2

Note that $Q'(y) = \frac{1}{f(Q(y))}$ and $Q''(y) = -\frac{f'(Q(y))}{f^3(Q(y))}$. Let $h(y) = (y(1-y))^\nu$. Moreover,

$$(h(y)\tilde{N}_n(y))' = h'(y)\tilde{N}_n(y) + h(y)\tilde{N}'_n(y).$$

Thus, in view of (21), we need $h'(y)$ and $h(y)Q'(y)$ to be uniformly bounded, which is achieved by $\nu > 1$ and $\nu > \gamma$, respectively. Moreover,

$$(h(y)\tilde{N}_n(y))'' = h''(y)\tilde{N}_n(y) + 2h'(y)\tilde{N}'_n(y) + h(y)\tilde{N}''_n(y).$$

Thus, in view of (21) and (22), we need $h''(y)$, $h'(y)Q'(y)$, $h(y)(Q'(y))^2$ and $h(y)Q''(y)$ to be uniformly bounded. The first claim is achieved by $\nu > 2$, the second by $\nu - 1 > \gamma$ and the third by $\nu > 2\gamma$. As for the fourth, we have

$$h(y)|Q''(y)| = \left(\frac{|f'(Q(y))|}{f^2(Q(y))}(y(1-y))\right)\left(\frac{(y(1-y))^{\nu-1}}{f(Q(y))}\right).$$

Since (CsR3(i)) and (CsR3(ii)) are fulfilled, the first part is $O(1)$ uniformly in $y \in (0,1)$. The second part is $O(1)$ since $\nu - 1 > 2\gamma - 1 > \gamma$ (recall that $\gamma > 1$).

Let $g(y) = h(y)\tilde{N}_n(y)$. Thus, proceeding exactly as in Lemma 3.5, we obtain $E[\max_{y \in (0,1)} \times |g'(y)|^\alpha] = O(n^{-\alpha/2})$. Thus, $\sup_{y \in (0,1)} |g'(y)| = o_{a.s}(l_q(n)/\sqrt{n})$ and

$$\sup_{|u-v| \le d_n, u,v \in (0,1)} |g(v) - g(u)| \le d_n \sup_{y \in (0,1)} |g'(y)| = o_{a.s}(l_q(n)/\sqrt{n}).$$

Further, by differentiability of $h$,

$$\sup_{|u-v| \le d_n, u,v \in (0,1)} |h(u) - h(v)\|\tilde{N}_n(u)| \le d_n \sup_{y \in (0,1)} |\tilde{N}_n(y)| = o_{a.s}(l_q(n)/\sqrt{n}).$$

Consequently,

$$\sup_{|u-v| \le d_n, u,v \in (0,1)} h(v)|\tilde{N}_n(v) - \tilde{N}_n(u)| = o_{a.s}(l_q(n)/\sqrt{n}).$$

Therefore, taking into account the martingale part and $d_n = \lambda_n$, we have

$$\sup_{|u-v| \le \lambda_n, u,v \in (0,1)} h(v)|(E_n(v) - v) - (E_n(u) - u)| = O_{a.s}\left(\frac{\sqrt{d_n \log n}}{\sqrt{n}}\right). \tag{23}$$

Note that

$$u_n(y) = \sqrt{n}(E_n(U_n(y)) - U_n(y)) + O_{a.s}(n^{-1/2}).$$

Consequently, by (23) and (9), we have

$$\sup_{y \in (0,1)} h(y)|u_n(y) - \alpha_n(y)| = \sup_{y \in (0,1)} h(y)|\alpha_n(U_n(y)) - \alpha_n(y)| + O_{a.s}(n^{-1/2}) = O_{a.s}(b_n).$$

Let $(k-1)/n < y \le k/n$. Further, let $\delta_n = 2C^* n^{-1/2} (\log\log n)^{1/2}$, $C^* = C$, where $C$ comes from (9). As in [5], with $\theta = \theta_n(y)$ such that $|\theta - y| \le n^{-1/2} u_n(y) = O_{a.s}(\delta_n)$,

$$\sup_{y \in (\delta_n, 1 - \delta_n)} h(y)|f(Q(y))q_n(y) - u_n(y)|$$

$$\le n^{-1/2} u_n^2(y) \left(\frac{f'(Q(\theta))}{f^2(Q(\theta))} \theta(1 - \theta)\right) \frac{f(Q(y))}{f(Q(\theta))} \frac{(y(1-y))^\nu}{\theta(1-\theta)}.$$

In view of Lemma 1 of [5], we have

$$\frac{f(Q(y))}{f(Q(\theta))} \le \left\{ \frac{y \vee \theta}{y \wedge \theta} \frac{1 - y \wedge \theta}{1 - y \vee \theta} \right\}^\gamma. \tag{24}$$

Further, if $y \ge 2\delta_n$, then

$$\frac{y}{\theta} = \frac{y - \theta}{\theta} + 1 \le \frac{\delta_n}{y - n^{-1/2} u_n(y)} + 1 \le 2. \tag{25}$$

Consequently, by (CsR3(i)), (CsR3(ii)), [24] and [25], we have

$$\sup_{y \in (\delta_n, 1-\delta_n)} h(y)|f(Q(y))q_n(y) - u_n(y)| = O_{a.s}(n^{-1/2} \log \log n).$$

If $1 \geq U_{k:n} \geq y$, then

$$
\begin{aligned}
\sup_{y \in (0, \delta_n)} h(y)|f(Q(y))q_n(y)| &= \sup_{y \in (0, \delta_n)} n^{1/2} h(y) \int_y^{U_{k:n}} \frac{f(Q(y))}{f(Q(u))} \, \mathrm{d}u \\
&\leq C \sup_{y \in (0, \delta_n)} n^{1/2} h(y) \int_y^{U_{k:n}} \left(\frac{u}{y}\right)^{\gamma_1} \mathrm{d}u \\
&\leq C \sup_{y \in (0, \delta_n)} n^{1/2} U_{k:n}^{\gamma_1+1} h(y) y^{-\gamma_1} \\
&= O_{a.s}(n^{-1/2} \log \log n).
\end{aligned}
$$

Further, if $U_{k:n} \leq y$, then for $\gamma_1 > 1$,

$$\sup_{y \in (0, \delta_n)} h(y)|f(Q(y))q_n(y)| \leq C n^{1/2} \delta_n^{\nu+\gamma_1} U_{k:n}^{-(\gamma_1-1)} \tag{26}$$

and for $\gamma_1 = 1$,

$$\sup_{y \in (0, \delta_n)} h(y)|f(Q(y))q_n(y)| \leq C n^{1/2} \delta_n^{\nu+1} \log(\delta_n/U_{k:n}). \tag{27}$$

Now,

$$P(U_{1:n} \leq n^{-2}(\log n)^{-3/2}) \leq \sum_{i=1}^n P(U_i \leq n^{-2}(\log n)^{-3/2}) \leq n^{-1}(\log n)^{-3/2}.$$

Consequently, via the Borel–Cantelli lemma, $U_{k:n} = o_{a.s}(n^{-2}(\log n)^{-3/2})$. Therefore, the bound in [27] is $O_{a.s.}(n^{-1/2} \log \log n)$ and the bound in [26] is of the same order, provided that $\nu > 3\gamma_1 - 2$. The upper tail is treated similarly. Consequently, the result follows. $\square$

# 5. Remarks

**Remark 5.1.** Assume that $c_k \sim k^{-\tau}(\log k)^{-3/2}$. Then,

$$\sum_{k=i}^\infty c_k^2 = \sum_{k=i}^\infty k^{-2\tau}(\log k)^{-3} = O(i^{-2\tau+1}(\log i)^{-3})$$

and [6] is fulfilled if $\tau > 5/2$.

If $X_1$ has all moments finite, then in [16], the summability assumption $\sum_{k=0}^\infty |c_k| < \infty$ provides the bound $O(n^{-3/4}(\log n)^{1/2+\eta})$, $\eta > 0$.

***Remark 5.2.*** It should be pointed out that the process $K(x,t)$ in Proposition 3.1 approximates $N_n$, not $\beta_n$. We use this as a tool to achieve ULIL only. Using the blocking technique applied to $\beta_n$ directly, we may establish Proposition 3.1 with $N_{[nt]}$ replaced by $([nt])^{1/2}\beta_{[nt]}$ and a different Gaussian process. However, it requires stronger assumptions than (6) on the coefficients $c_k$. In particular, $\tau > \frac{9}{2} + \frac{4}{\alpha}$. Without exploiting differentiability, we must compare indicators applied to the original sequences $X_i$ with the indicators applied to the truncated one (see, e.g., [3]).

***Remark 5.3.*** The main part of this paper is devoted to the law of the iterated logarithm for the empirical process. One may ask whether, by imposing (6), we could obtain uniform law of the iterated logarithm for empirical processes via strong mixing. From [15], we know that if $\mathrm{E}|\epsilon_1|^\alpha < \infty$, $c_k = O(k^{-\tau})$ and $\tau > 1 + \frac{1}{\alpha} + \max\{1, \alpha^{-1}\}$, then $\alpha(n) = O(n^{-\lambda})$, where $\lambda = (\tau\alpha - \max\{\alpha,1\})/(1+\alpha) - 1 > 0$ and $\alpha(n)$ is the strong mixing coefficient. In particular, if $\alpha \geq 2$, then this condition yields $\tau > 2 + \frac{1}{\alpha}$ and $\alpha(n) = O(n^{-(\tau\alpha-1-2\alpha)/(1+\alpha)})$. In view of Rio [14], to obtain the LIL for partial sums of bounded strongly mixing random variables, we need $\alpha(n)$ to be summable, which would require $\tau > 3 + \frac{2}{\alpha}$. Thus, we need stronger conditions than (6). Also, Rio's result would provide the LIL for $\beta_n(x)$ at *fixed* $x$.

# Appendix

Here, we reprove the strong approximation of Proposition 3.1. The proof is along the lines of [3], thus we present only the major steps.

We keep the notation from Section 3.1. Recall that $\|X\|_\alpha = (\mathrm{E}|X|^\alpha)^{1/\alpha}$. For $\mathbf{x} = (x_1, \ldots, x_d) \in \mathbb{R}^d$, let $\|\mathbf{x}\|$ be the Euclidean norm and define the following random vectors in $\mathbb{R}^d$:

$$\xi_i = (Y_i(x_1), \ldots, Y_i(x_d)), \qquad \hat{\xi}_i = (\hat{Y}_i(x_1), \ldots, \hat{Y}_i(x_d)).$$

Denote the covariance matrix of $\xi_i$ by $\Gamma_d = \Gamma_d(\mathbf{x}) = (\Gamma(x_i, x_j), 1 \leq i, j \leq d)$. We shall use the following implication, which is valid for any vectors $\xi_i$, $\eta_i$ in $\mathbb{R}^d$:

$$\mathrm{E}\|\xi_i - \eta_i\| \leq A \quad \text{implies} \quad |\mathrm{E}\exp(i\langle \mathbf{u}, \xi_i \rangle) - \mathrm{E}\exp(i\langle \mathbf{u}, \eta_i \rangle)| \leq \|\mathbf{u}\|A. \tag{28}$$

Also, we will use the following bound: for $\mathbf{x}_i \in \mathbb{R}^d$, $i = 1, \ldots, n$, we have $\|\sum_{i=1}^n \mathbf{x}_i\| \leq \sum_{i=1}^n \|\mathbf{x}_i\|$.

**Lemma A.1.** *Under the conditions of Proposition 3.1, for all* $\mathbf{u} \in \mathbb{R}^d$,

$$\left| \mathrm{E}\exp\left( i\left\langle \mathbf{u}, n^{-1/2}\sum_{i=1}^n \xi_i \right\rangle \right) - \mathrm{E}\exp\left( i\left\langle \mathbf{u}, n^{-1/2}\sum_{i=1}^n \hat{\xi}_i \right\rangle \right) \right| \leq \|\mathbf{u}\|d^{1/2}n^{-1/2}.$$

**Proof.** From (6) and (17), we have $\|Y_i(x) - \hat{Y}_i(x)\|_2 \leq Ci^{-1}(\log i)^{-3/2}$. Thus,

$$\mathrm{E}\|\xi_i - \hat{\xi}_i\|^2 \leq Cdi^{-2}(\log i)^{-3} \tag{29}$$

and (29) implies that

$$\mathrm{E}\|\xi_i - \hat{\xi}_i\| \le C d^{1/2} i^{-1} (\log i)^{-3/2}.$$

Therefore,

$$\mathrm{E}\left\|\sum_{i=1}^{n}(\xi_i - \hat{\xi}_i)\right\| \le C d^{1/2}.$$

This, together with (28), implies the result. □

**Proposition A.2.** *Under the conditions of Theorem 3.1, for all* $\mathbf{u} \in \mathbb{R}^d$,

$$\left|\mathrm{E}\exp\left(i\left\langle\mathbf{u}, n^{-1/2}\sum_{i=1}^{n}\xi_i\right\rangle\right) - \exp(-\tfrac{1}{2}\langle\mathbf{u}, \Gamma_d\mathbf{u}\rangle)\right|$$

$$\le C d\|\mathbf{u}\|n^{-\delta_1} + C\|\mathbf{u}\|^2(n^{-1/4}(\log n)^{d/2}\exp(Cd) + d^{d/2}\exp(-Cn^{1/2}))$$

*for some* $\delta_1 > 0$.

**Proof.** Divide the interval $[1, n]$ into consecutive long and short blocks $I_1, J_1, I_2, J_2, \ldots$ with lengths $card(I_k) = [n^{\rho^*}]$, $card(J_k) = [n^\rho]$, $1 \le k \le M$, $\rho < \rho^* < \frac{1}{2}$. Then, $M = M_n \sim n^{1-\rho^*}$ (the last block is possibly incomplete). Define $U_k = \sum_{i \in I_k}\xi_i$, $\hat{U}_k = \sum_{i \in I_k}\hat{\xi}_i$ and $\hat{V}_k = \sum_{i \in J_k}\hat{\xi}_i$. By (29), we have

$$\mathrm{E}\|U_k - \hat{U}_k\| \le \sum_{i \in I_k}\mathrm{E}\|\xi_i - \hat{\xi}_i\| \le C d^{1/2}\sum_{i \in I_k} i^{-1}(\log i)^{-3/2}. \tag{30}$$

Note that $\hat{U}_k$, $k = 1, \ldots, M$, are independent. Thus, we can construct independent random vectors $\tilde{U}_1, \ldots, \tilde{U}_M$ such that $(\tilde{U}_k, \hat{U}_k) \overset{\mathrm{d}}{=} (U_k, \hat{U}_k)$, $k = 1, \ldots, M$. By (30), we obtain

$$\mathrm{E}\left\|\sum_{k=1}^{M}(\tilde{U}_k - \hat{U}_k)\right\| \le \sum_{k=1}^{M}\mathrm{E}\|U_k - \hat{U}_k\| \le C d^{1/2}. \tag{31}$$

Further, $\mathrm{E}\|\hat{V}_k\| \le dn^{\rho/2}$, which implies, by independence, that

$$\mathrm{E}\left\|\sum_{k=1}^{M}\hat{V}_k\right\|^2 \le \sum_{k=1}^{M}\mathrm{E}\|\hat{V}_k\|^2 \le C d M n^\rho = d n^{\rho+1-\rho^*}. \tag{32}$$

We have

$$\sum_{i=1}^{n}\hat{\xi}_i = \sum_{k=1}^{M}(\hat{U}_k - \tilde{U}_k) + \sum_{k=1}^{M}\tilde{U}_k + \sum_{k=1}^{M}\hat{V}_k.$$

Thus, by (31) and (32)

$$\mathrm{E}\left\|\sum_{i=1}^{n}\xi_i - \sum_{k=1}^{M}\tilde{U}_k\right\| \le C(d^{1/2} + dn^{(\rho+1-\rho*)/2}). \tag{33}$$

Also, by Lemma 2.9 of [3] and the same argument as at the end of Lemma 2.10 in the same paper, we obtain

$$\left|\mathrm{E}\exp\left(i\left\langle \mathbf{u}, n^{-1/2}\sum_{k=1}^{M}\tilde{U}_k\right\rangle\right) - \exp\left(-\tfrac{1}{2}\langle\mathbf{u},\Gamma_d\mathbf{u}\rangle\right)\right|$$

$$\le C\|\mathbf{u}\|^2(n^{-(1-\rho^*)/2}(\log n)^{d/2}\exp(Cd) + d^{d/2}\exp(-Cn^{1-\rho^*})). \tag{34}$$

Consequently, by (33), (28) and (34), the result follows.                                    □

## A.1. Approximation

Let $\varepsilon \in (0, 1/4)$. Further, let $t_k = \exp(k^{1-\varepsilon})$, $p_k = 2[t_k^{\rho}]$, $d_k = k^{1/2}$ and $x_i = (i-1)/d_k$, $i \in \mathbb{Z}$. Define $M_k = t_{k+1} - t_k - p_k$ so that $M_k \sim Ck^{-\varepsilon}\exp(k^{1-\varepsilon})$ and $\log M_k \sim Ck^{1-\varepsilon}$. Define random vectors in $\mathbb{R}^{d_k}$, $\eta_k = (\eta_{k1}, \ldots, \eta_{kd_k})$, where

$$\eta_{ki} = R(x_i, t_{k+1}) - R(x_i, t_k + p_k) = \sum_{j=t_k+p_k+1}^{t_{k+1}} Y_j(x_i).$$

Also, define $\hat{\eta}_{ki} = \sum_{j=t_k+p_k+1}^{t_{k+1}} \hat{Y}_j(x_i)$.

By Proposition A.2,

$$|\mathrm{E}\exp(i\langle\mathbf{u}, M_k^{-1/2}\eta_k\rangle) - \exp(-\tfrac{1}{2}\langle\mathbf{u},\Gamma_{d_k}\mathbf{u}\rangle)|$$

$$\le C\|\mathbf{u}\|^2(d_k M_k^{-\rho_1} + M_k^{-(1-\rho^*)/2}(\log M_k)^{d_k/2}\exp(Cd_k) + d_k^{d_k/2}\exp(-CM_k^{1-\rho^*}))$$

$$\le C\exp(-Ck^{1-\varepsilon})\|\mathbf{u}\|^2.$$

Let $\psi_{PL}$ be the Prokhorov–Lévy distance. Choose $T := \exp(k^{\varepsilon})$; for sufficiently large $k$, $T > 10^8 d_k$. Then,

$$\psi_{PL}(M_k^{-1/2}\eta_k, N(\mathbf{0}, \Gamma_{d_k}))$$

$$\le \frac{16 d_k}{T}\log T + P(N(\mathbf{0}, \Gamma_{d_k}) > T/2)$$

$$\quad + T^{d_k}\int_{\|\mathbf{u}\|\le T}\left|\mathrm{E}\exp(i\langle\mathbf{u}, M_k^{-1/2}\eta_k\rangle) - \exp\left(-\frac{1}{2}\langle\mathbf{u},\Gamma_{d_k}\mathbf{u}\rangle\right)\right|\mathrm{d}\mathbf{u}$$

$$\le \frac{16 d_k}{T}\log T + P(N(\mathbf{0}, \Gamma_{d_k}) > T/2) + CT^{d_k}\exp(-k^{1-\varepsilon})\int_{\|\mathbf{u}\|\le T}\|\mathbf{u}\|^2\,\mathrm{d}u$$

$$\le C \exp(-Ck^\varepsilon)$$

since $\varepsilon < 1/4$.

Since $\hat{\eta}_k$, $k = 1, \ldots, M$, are independent, we can define independent random vectors $\zeta_1, \ldots, \zeta_M$ such that $M_k^{-1/2} \zeta_k \sim N(\mathbf{0}, \Gamma_{d_k})$ and

$$P(\|M_k^{-1/2} \hat{\eta}_k - M_k^{-1/2} \zeta_k\| > C \exp(-k^\varepsilon)) \le C \exp(-k^\varepsilon).$$

This yields

$$P(\|M_k^{-1/2} \eta_k - M_k^{-1/2} \zeta_k\| > C \exp(-k^\varepsilon)) \le C \exp(-k^\varepsilon)$$

and then, by the Borel–Cantelli lemma,

$$M_k^{-1/2}(\eta_k - \zeta_k) \le \exp(-k^\varepsilon) \qquad \text{almost surely.} \tag{35}$$

For $k \ge 1$, define random vectors $Z_k$ in $\mathbb{R}^{d_k}$ by

$$Z_k = (R(x_i, t_{k+1}) - R(x_i, t_k), i = 1, \ldots, d_k).$$

Since $Z_k - \eta_k$ is the sum of $p_k$ random vectors in $\mathbb{R}^{d_k}$ with coordinates bounded by 1 and since $M_k \le t_k$, $d_k^{1/2} p_k^{1/2} \le k t_k^{\rho/2} \le \exp(ck^{1-\varepsilon})$, $\varepsilon < 1/2$, we have, by (35),

$$\|Z_k - \zeta_k\| \le \|Z_k - \eta_k\| + \|\eta_k - \zeta_k\|$$
$$\le C(d_k^{1/2} p_k^{1/2} + t_k^{1/2} \exp(-k^\varepsilon)) \le C t_k^{1/2} \exp(-Ck^\varepsilon). \tag{36}$$

Thus, the skeleton process $\{Z_k, k \ge 1\}$ can be approximated by the sequence $\{\zeta_k, k \ge 1\}$. The latter can be extended to a centered Gaussian process $\{K(x, t), x \in \mathbb{R}, t \ge 0\}$ with covariance $t \wedge t' \Gamma_{d_k}(s, x')$ such that

$$\zeta_{ki} = K(x_i, t_{k+1}) - K(x_i, t_k + p_k), \qquad i = 1, \ldots, d_k.$$

Define $Y_k = (K(x_i, t_{k+1}) - K(x_i, t_k), i = 1, \ldots, d_k)$. Then,

$$\|\zeta_k - Y_k\| \le C t_k^{1/2} \exp(-Ck^\varepsilon) \qquad \text{almost surely}$$

by the last inequality on page 807 of [3]. Thus, by (36),

$$\|Z_k - Y_k\| \le C t_k^{1/2} \exp(-Ck^\varepsilon) \qquad \text{almost surely.}$$

## A.2. Oscillations

**Lemma A.3.** *Under the conditions of Theorem 3.1, for $n \ge 1$, $\lambda \ge n^{1/2}$ and any $-\infty < a < b < \infty$, we have*

$$P\left(\sup_{a \le x \le x' \le b, 0 \le k \le n} \sum_{i=1}^{k} Y_i(x, x') \ge \lambda\right) \le C \exp(-C\lambda^2/(n(b-a)^2)) + \frac{C}{n^\eta},$$

*where $\eta > 0$.*

**Proof.** Define

$$M_{u,v} = \max_{0 \leq i \leq 2^u, 0 \leq j \leq 2^v} \left| \sum_{k=nj2^{-v}}^{n(j+1)2^{-v}} Y_k(bi2^{-u}, b(i+1)2^{-u}) \right|.$$

Then, for any integer $L \geq 1$,

$$\left| \sum_{i=1}^{k} Y_i(0, x) \right| \leq \sum_{1 \leq u, v \leq L} M_{u,v} + \frac{2n}{2^L}.$$

Choose

$$0 < \varepsilon < (\alpha - 2)/(2C_6) \tag{37}$$

for some constant $C_6$ to be specified later. Since $n^{1/2} \leq \lambda$, $u, v \leq L$, we have, by Lemma 3.4,

$$P(M_{u,v} \geq \lambda 2^{-\varepsilon(u+v)})$$
$$\leq 2^{u+v} C \{ (\lambda 2^{-\varepsilon(u+v)})^{-\alpha} + \exp(-C\lambda^2 2^{-2\varepsilon(u+v)}/(n2^{-v}(b2^{-u})^2)) $$
$$+ \exp(-C\lambda 2^{-\varepsilon(u+v)}/(n2^{-v})^\rho) \}$$
$$\leq 2^{u+v} C \{ n^{-1/2\alpha} 2^{2\varepsilon L\alpha} + \exp(-C\lambda^2 2^{u(1-2\varepsilon)} 2^{v(1-2\varepsilon)}/(nb^2)) $$
$$+ (n^{(1/2-\rho)} 2^{-\varepsilon(u+v)} 2^{v\rho})^{-r} \}$$

for any $r > 0$. Choose $r = \frac{\alpha}{2(1/2-\rho)}$. Since $u + v \leq 2L$, we have

$$P(M_{u,v} \geq \lambda 2^{-\varepsilon(u+v)})$$
$$\leq 2^{u+v} C \{ n^{-\alpha/2} 2^{2C\varepsilon L\alpha} + \exp(-C\lambda^2 2^{u(1-2\varepsilon)} 2^{v(1-2\varepsilon)}/(nb^2)) \}.$$

For

$$A = \{ M_{u,v} \geq \lambda 2^{-\epsilon(u+v)} \text{ for some } 1 \leq u, v \leq L \},$$

we have

$$P(A) \leq \sum_{1 \leq u, v \leq L} 2^{u+v} \exp(-C\lambda^2 2^{u(1-2\varepsilon)} 2^{v(1-2\varepsilon)}/(nb^2)) + 2^{2L+2} 2^{CL\varepsilon} n^{-\alpha/2}$$

$$\leq \left( \sum_{1 \leq u \leq L} \exp(u - C\lambda^2 2^{u(1-2\varepsilon)}/(nb^2)) \right)^2 + 2^{L(2+C\varepsilon)+2} n^{-\alpha/2}$$

$$\leq C \exp(-C\lambda^2/(nb^2)) + C 2^{L(2+C_6\varepsilon)} n^{-\alpha/2}.$$

If we choose $L$ so that $n^{1/2} < 2^L < 2n^{1/2}$, then, by (37), the second part in the latter expression is bounded above by $Cn^{-\eta}$, for some $\eta > 0$. Further, on $A^c$, we have $\sum_{1 \le u,v \le L} M_{u,v} \le C_7\lambda$ and

$$\sum_{1 \le u,v \le L} M_{u,v} + \frac{2n}{2^L} \le C_7\lambda + 4\lambda = C_8\lambda.$$

Therefore, the result follows. $\qquad\square$

Let

$$\tilde{R}_{i,k} = \sup_{s \in [x_i, x_{i+1})} \sup_{t \in [t_k, t_{k+1}]} |R(s,t) - R(x_i, t_k)|.$$

Using Lemma A.3, we obtain the following result.

**Lemma A.4.** *Under the conditions of Theorem 3.1, for some $\varepsilon^* > 0$,*

$$\max_{1 \le i \le d_k} |\tilde{R}_{i,k}| \le Ct_k^{1/2}(\log k)^{-\varepsilon^*} \qquad \text{almost surely.}$$

## A.3. Conclusion of the proof

For $s \in [x_i, x_{i+1})$ and $t \in [t_k, t_{k+1}]$, write

$$|R(s,t) - K(s,t)| \le |R(s,t) - R(x_i, t_k)| + |K(s,t) - K(x_i, t_k)| + |R(x_i, t_k) - K(x_i, t_k)|.$$

Both the first and the second part are bounded almost surely by $Ct_k^{1/2}(\log k)^{-\varepsilon^*}$, the first by Lemma A.4 and the second by (2.55) in [3]. The second part is bounded in the same way as (2.56) in [3].

# Acknowledgements

The comments of Neville Weber and both referees are greatly appreciated. Part of this work was done during my stay at Carleton University. I am thankful to Professors Barbara Szyszkowicz and Miklos Csörgő for their support and helpful remarks.